\documentclass[12pt]{article}
\usepackage{amssymb, amsmath, amsfonts}
\font\tenfrak=eufm10
\font\sevenfrak=eufm7
\font\fivefrak=eufm5
        \newfam\frakfam 
                \textfont\frakfam=\tenfrak
\font\tendl=msbm10  scaled \magstep1%double line
\font\sevendl=msbm7 scaled \magstep1
\font\fivedl=msbm5 scaled \magstep1
        \newfam\dlfam  % \dl is double line
                \textfont\dlfam=\tendl
\scriptfont\dlfam=\sevendl
        \scriptscriptfont\dlfam=\fivedl
\scriptfont\frakfam=\sevenfrak
        \scriptscriptfont\frakfam=\fivefrak
\long\def\nodo#1{}
\def\genfd{{\bf k}}
\def\id{{\rm id}}

\newcounter{rsec}
\newcounter{rssec}
\newcounter{point}

\newcommand{\ppt}{ \addtocounter{point}{1}
{\bf \arabic{section}.\arabic{point} }}
\newcommand{\ppta}{
{\bf \arabic{section}.\arabic{point}a }}
\newcommand{\pptb}{
{\bf \arabic{section}.\arabic{point}b }}
\newcommand{\pptc}{
{\bf \arabic{section}.\arabic{point}c }}
\newcommand{\pptd}{
{\bf \arabic{section}.\arabic{point}d }}

\def\lnameuse#1{\expandafter\csname lpt#1 \endcsname}
\def\myedef#1{\expandafter\edef\csname #1\endcsname}%%I want edef instead
\def\muse#1{\expandafter\csname #1\endcsname}%%avoids @

\begin{document}
\begin{center}
{\large \bf Quantum heaps, cops and heapy categories}\\
\vskip .03in
\footnotesize{
Zoran \v{S}koda, {\tt zskoda@irb.hr} \\
Institute Rudjer Bo\v{s}kovi\'{c}, P.O.Box 180,\\
HR-10002 Zagreb,
Croatia
}
\end{center}
%Subj-class: Quantum Algebra; Rings and Algebras
%MSC-class: 20N10, 16A24

%\begin{abstract}
A heap is a structure with a ternary operation which is intuitively
a group with forgotten unit element. Quantum heaps are associative
algebras with a ternary cooperation which are to
the Hopf algebras what heaps
are to groups, and, in particular, the category of copointed
quantum heaps is isomorphic to the category of Hopf algebras.
There is an intermediate structure of a cop in monoidal
category which is in the case of vector spaces
to a quantum heap about what is a coalgebra to a Hopf algebra.
The representations of Hopf algebras make a rigid monoidal
category. Similarly the representations of quantum heaps
make a kind of category with ternary products, which we call
a heapy category.

%\end{abstract}

%\addtocounter{chapter}{1}
\section{Classical Background: heaps}
 \addtocounter{rsec}{1}\setcounter{point}{0}

%\scdraft{qheapthes.tex}
A reference for this section is~\cite{bergman:corings}.
\nodo{
 \ppt {\it Motivation.}
Heap is an algebraic structure which is basically
equivalent to group when one forgets which element is unit.

A similar idea is encountered in the concept of principal
homogeneous set of a group, that is a set on which a given
group acts freely and transitively. The group
in the heap theory will appear as a natural automorphism
group.
}

\ppt {\it Before forgetting: group as heap.}
Let $G$ be a group. Then the ternary
operation $t : G \times G \times G \rightarrow G$ given by
\begin{equation}\label{eq:t-op} t(a,b,c) = ab^{-1}c, \end{equation}
satisfies the following relations:
\begin{equation}\label{eq:def-heap}
\begin{array}{l} t(b,b,c) = c = t(c,b,b)\\
        t(a,b,t(c,d,e)) = t(t(a,b,c),d,e) \end{array}\end{equation}
\nodo{
We interpret that operation as shifting $a$ by the (right)
translation in the group which translates $b$ into $c$.
}
 A heap $(H,t)$ is a pair of the nonempty set $H$
and a ternary operation $t : H \times H \times H$ satisfying
relation~(\ref{eq:def-heap}). A morphism $f : (H,t)\to (H',t')$
of heaps is a set map $f : H\to H'$
satisfying $t'\circ (f\times f\times f) = f \circ t$.

Thus every group has its canonical heap, what defines a
faithful functor Heap : \underline{Groups} $\rightarrow$
\underline{Heaps}.

\ppt The automorphism group of a heap $(H,t)$ denoted by ${\rm Aut}H$
is the subgroup of the symmetric group of $H$ consisting of the maps
of the form $t(\cdot,a,b): H \rightarrow H$
where $a,b \in H$, and $\cdot$ is a place-holder.
The composition (group operation) satisfies
\[ t(\cdot,c,d) \cdot_{{\rm Aut}H} t(\cdot, a,b) =
 t(t(\cdot,c,d),a,b) = t(\cdot,c,t(d,a,b)). \]
The rightmost equality implies that the result of the composition
is in ${\rm Aut}H$.
%If that is (what is most sensible) taken as
%the definition then the identity above is axiom for the action.
The inverse
of $t(\cdot,a,b)$ is $t(\cdot,b,a)$ by~(\ref{eq:def-heap})
and the unit is $t(\cdot, x,x)$ (independent of $x \in H$).

\ppt The following are equivalent

(i) bijections $t(\cdot,a,b)$ and $t(\cdot,a',b')$ are
the same maps,

(ii) $t(a,a',b') = b$,

(iii) $t(b,b',a') = a$.

Proof. (ii) follows from (i) and $t(a,a,b) = b$.

(iii) follows from (ii) by applying $t(\cdot,b',a')$ on the
right. Similarly (ii) follows from (iii).

(i) follows from (ii) by the calculation:
\[ t(x,a',b') = t(t(x,a,a),a',b')= t(x,a,t(a,a',b'))
 = t(x,a,b).\mbox{ }\]

\ppt The reader should conclude noticing that
the defining action of ${\rm Aut}H$ is transitive
(by $t(a,a,b) = b$) and  free (if $t(a,b,c) = a$ then
by lastpt $t(x,b,c) = x$ for each $x$, in particular
$t(b,b,c) = b$ but it also $t(b,b,c) = c$ by~(\ref{eq:def-heap}).

\ppt If we started with a group $G$ then we can recover it
{\it up to isomorphism}
from the corresponding heap (as the automorphism group of
the heap).
 Indeed, then the $t(\cdot,e,a)$ is the multiplication by
$a$. A byproduct of this construction is that we now know
that all the other possible identities for the group-induced
ternary operation~(\ref{eq:t-op}) follow from~(\ref{eq:def-heap}).

\ppt Similarly every heap is isomorphic (in the category of heaps,
where morphisms are defined as it is usual for algebraic
structures) to the heap of operation~(\ref{eq:t-op})
on its automorphism group. However the isomorphism is not
natural but one needs to specify which element will be unity.
In other words we have the natural isomorphism (not only equivalence)
of the category of the groups with the category of {\bf pointed heaps},
that is heaps with a nullary operation $\star$ and morphisms
respecting also this operation.

The isomorphism in question is $H \ni a \mapsto t(\cdot,\star,a)
\in \iota(H)= {\rm Heap}({\rm Aut}H)$.
\[\begin{array}{l} t(\cdot,\star,a) [t(\cdot,\star,b)]^{-1} t(\cdot,\star,c)
= t(\cdot,\star,a) t(\cdot,b,\star) t(\cdot,\star,c) \\ =
t(\cdot,\star,a)t(\cdot,b,c) = t(\cdot,\star,t(a,b,c)).
\end{array}\]

For the morphism of heaps $f : (H,t) \rightarrow (H',t')$
we define $\iota(f)(t(\cdot,\star,a))= t'(\cdot,\star,f(a))$
and $\iota$ becomes a covariant functor.

The identities for the ternary operation $t$ play important role in
universal algebra (theory of Mal'cev algebras and Mal'cev categories).

\section{Cops}

\ppt Let $({\cal C},\otimes,{\bf 1})$, or ${\cal C}$ for short,
be a strict monoidal category with unit object ${\bf 1}$.
A {\bf cop} $C$ in $({\cal C},\otimes, {\bf 1})$ is a pair
$(C,\tau)$, where $C$ is an object in ${\cal C}$
and $\tau : C \rightarrow C \otimes C \otimes C$
a morphism in ${\cal C}$ satisfying the law
\begin{equation}\label{eq:cop}
(\id \otimes \id \otimes \tau) \circ \tau =
  (\tau \otimes \id \otimes \id) \circ \tau.
\end{equation}
Let $({\cal C},\otimes,{\bf 1},\sigma)$ be a strict symmetric monoidal
category and $C$ a monoid (=algebra) object in that category, i.e.
$C$ is equipped with a product $\mu : C \otimes C \rightarrow C$
and a unit morphism $\eta : {\bf 1} \rightarrow C$ satisfying the
standard axioms.
Then an {\bf opposite monoid} $C_{{\rm op}}$ is the same object $C$
equipped with product $\sigma_{C,C} \circ \mu$ and
with the same unit map $\eta$.
A {\bf symmetric cop monoid} $C$ in
a strict symmetric monoidal category $({\cal C},\otimes,{\bf 1},\sigma)$
is a {\it monoid} object $C$ with a {\it morphism of monoids}
$\tau : C \rightarrow C \otimes C_{{\rm op}} \otimes C$
satisfying the law~(\ref{eq:cop}). Here the tensor product has
the usual tensor product
structure of a monoid in a strict symmetric monoidal category
(for two monoids $A$ and $B$ one takes $(\mu \otimes \mu) \circ
(\id \otimes \sigma_{B,A} \otimes \id)$ as a product on $A \otimes B$).
%for three monoids one needs to use bracketing, but any choice will do).

\ppt A {\bf counit} of a cop $C$ in ${\cal C}$ is a morphism
$\epsilon : C \rightarrow {\bf 1}$ such that
\begin{equation}\label{eq:counit}
 (\id \otimes \epsilon \otimes \epsilon) \circ \tau = \id = (\epsilon
\otimes \epsilon \otimes \id) \circ \tau ,
\end{equation}
where the identification morphism
${\bf 1} \otimes {\bf 1} \otimes C \equiv C
\equiv C \otimes {\bf 1} \otimes {\bf 1}$ is used.
%(with the same
%bracketings as the ones implicit in equation~(\ref{eq:counit})).

Our interest is in the cops in the (symmetric) monoidal
category of vector spaces $\underline{Vec}$ or supervector spaces
$\underline{SVec}$ over some fixed field.

\ppt A {\bf coheap monoid} in a symmetric monoidal category is a
symmetric cop monoid such that
$(\id \otimes \mu)\circ \tau = (\mu\otimes \id)\circ\tau = \id$,
where the identification $C\otimes {\bf 1} = C = {\bf 1}\otimes C$
is implicitly used.

\ppta A {\bf character} of a monoid $(C,\mu,\eta)$
in a strict monoidal category
is a morphism $\epsilon : C \rightarrow {\bf 1}$ such that
$\epsilon \circ \eta = \id_{\bf 1}$ and
$(\epsilon \otimes \epsilon) = \epsilon \circ\mu$.
\nodo{
and the following diagram commutes
\[ \begin{array}{ccc}
C \otimes C &
\stackrel{\epsilon \otimes \epsilon}{\rightarrow} &
{\bf 1}\otimes {\bf 1}\\
\,\,\,\,\downarrow \mu && \,\,\,\,\,\downarrow {\rm =} \\
C & \stackrel\epsilon\rightarrow & {\bf 1}\,\,\,.
\end{array}\]
}

\pptb A {\bf character} of a symmetric cop monoid $C$
is any character of $(C,\eta,\mu)$ in ${\cal C}$.

\pptc {\it Proposition.} A {\bf character} of a coheap monoid
is a automatically a counit of the underlying cop.

Proof is straightforward:
$$\begin{array}{c}
(\id \otimes \epsilon\otimes \epsilon)\tau =
(\id\otimes (\epsilon\circ\mu))\tau =
(\id\otimes\epsilon)(\id\otimes\mu)\tau =
(\id\otimes\epsilon)(\id\otimes\eta) = \id,
\\
(\epsilon\otimes\epsilon\otimes\id)\tau =
((\epsilon\circ\mu)\circ\id)\tau = (\id\otimes\epsilon)(\mu\otimes\id)\tau =
(\epsilon\otimes\id)(\eta\otimes\id) = \id,
\end{array}$$
where again obvious identifications are implicitly
used, e.g. $\id_{\bf 1} \otimes \id  \cong \id$.

\pptd A {\bf copointed cop} is a pair $(C,\epsilon)$
of a cop $C$ and a counit $\epsilon$ of $C$.
A {\bf copointed coheap monoid} is a coheap monoid with
a {\it character} $\epsilon$ of $C$. Warning: a copointed cop which is also
a coheap is not necessarily a copointed coheap, as the counit does not
need to be a character of a coheap. Clearly the above theory may be
modified for nonstrict monoidal categories.

\section{Quantum heaps}
  \addtocounter{rsec}{1}\setcounter{point}{0}

\indent
\ppt Heap is morally a group with forgotten unit.
Quantum heap is morally a Hopf algebra with forgotten counit.
We fix a ground field $\genfd$ throughout.

\ppt {\bf Quantum heap} is an associative
unital $\genfd$-algebra $(H,\mu,\eta)$
together with a ternary algebra cooperation
\[ \tau : H \rightarrow H \otimes H_{{\rm op}} \otimes H,\]
satisfying the following properties
\begin{equation}\begin{array}{l}
(\id  \otimes \id \otimes \tau) \tau
= (\tau \otimes \id \otimes \id) \tau \\
(\id \otimes \mu) \tau = \id \otimes 1_H \\
(\mu \otimes \id) \tau = 1_H \otimes \id
\end{array}\end{equation}
Moreover, $\tau$ is required to be algebra homomorphism from $H$ into
$H \otimes H_{{\rm op}} \otimes H$, where  $H_{{\rm op}}$ has the opposite
algebra structure and the tensor product has the usual algebra structure.
In other words, it is a coheap monoid in the symmetric monoidal category of
vector spaces.

We use heap analogue of the Sweedler notation:
\[ \tau(h) = \sum h^{(1)} \otimes h^{(2)} \otimes h^{(3)},\]
and because of the first of the above identities,
it is justified to extend it to any odd number $\geq 3$ factors, e.g.
\[ (\id \otimes \id \otimes \tau) \tau(h) =
 \sum h^{(1)} \otimes h^{(2)} \otimes h^{(3)}
\otimes h^{(4)} \otimes h^{(5)}. \]
In this paper heap-Sweedler notation has {\it upper} indices while the Sweedler
notation for coalgebras will have {\it lower} indices. Heap-Sweedler indices
extend to any any odd number $\geq 3$ of indices.
The requirement that $\tau$ is an algebra homomorphis from $H$ into
$H \otimes H_{{\rm op}} \otimes H$ is expressed in terms of
heap-Sweedler notation as
\[ \tau(hg) = \sum (hg)^{(1)} \otimes (hg)^{(2)}\otimes (hg)^{(3)}
= \sum h^{(1)}g^{(1)} \otimes g^{(2)}h^{(2)} \otimes h^{(3)} g^{(3)}
= \tau(h)\tau(g) \]

\ppt A morphism of quantum heaps is a homomorphism $\phi$ of unital algebras
such that $\tau(\phi(h)) = (\phi \otimes \phi \otimes \phi)\tau(h)$.
Quantum heaps make a category $\underline{QHeaps}$.

\ppt We define a covariant functor $QHeap$ from the category of Hopf algebras
$\underline{Hopf-Alg}$ to $\underline{QHeaps}$.
The underlying associative algebra of the object is the same,
and the quantum heap operation is given by
\[
 \tau(h) = \sum h_{(1)} \otimes Sh_{(2)} \otimes h_{(3)},
\]
This functor is identity on morphisms. However not all morphisms
of quantum heaps are in the image of this functor (for example
take a coordinate ring of $SL(n)$ and permute the rows of the matrix
of generators -- it will induce a morphism between the canonical and
the obvious ``permuted'' quantum heap structures).

Let us prove that this functor has the required codomain i.e.
indeed the output of functor $QHeap$ is in $\underline{QHeaps}$:
\[\begin{array}{l}
(\id \otimes \id \otimes \tau) \tau (h) =
\sum h_{(1)} \otimes Sh_{(2)} \otimes
( h_{(3)} \otimes Sh_{(4)} \otimes h_{(5)}) \\
\,\,\,= \sum ( h_{(1)} \otimes Sh_{(2)} \otimes
h_{(3)} ) \otimes Sh_{(4)} \otimes h_{(5)} =
(\tau \otimes \id \otimes \id) \tau(h)
\end{array}\]
\[\begin{array}{l}
(\id \otimes \mu)\tau(h) = \sum h_{(1)} \otimes Sh_{(2)} h_{(3)}
= h \otimes 1\\
(\mu \otimes \id)\tau(h) = \sum h_{(1)} Sh_{(2)} \otimes h_{(3)}
= 1 \otimes h
\end{array}\]
\[\begin{array}{lcl}
\tau(hg)&=& \sum (hg)_{(1)} \otimes S((hg)_{(2)}) \otimes (hg)_{(3)}\\
&=&  \sum h_{(1)}g_{(1)} \otimes (S(g_{(2)})\cdot_H S(h_{(2)}))
\otimes h_{(3)}g_{(3)}\\
&=&  \sum h_{(1)}g_{(1)} \otimes (S(h_{(2)})\cdot_{H_{\rm op}} S(g_{(2)}))
\otimes h_{(3)}g_{(3)} \\
& = & \tau(h)\tau(g)
\end{array}\]

\ppt Let now ${\bf A} = (A, \mu, \eta, \tau)$ be a quantum heap
and $\epsilon : A \rightarrow {\bf k}$ any
unital ${\bf k}$-algebra homomorphism.
We call such a pair $({\bf A},\epsilon)$ a {\bf copointed quantum heap}.
A morphism of copointed quantum heaps is a morphism
$\phi : {\bf A}\rightarrow {\bf A'}$
of quantum heaps such that $\epsilon' \circ \phi = \epsilon$.
Copointed quantum heaps make a category $\underline{Q\star Heaps}$
(notation from view of them as quantum ``pointed heaps'').

Now we define a functor
$\chi :\underline{Q\star Heaps}\rightarrow \underline{Hopf-Alg}$.

For a given heap $H$ thus define
\[ \Delta : H \rightarrow H \otimes {\bf k} \otimes H \cong H \otimes H
\,\,\,\,\,\mbox{ by }\,\,\,\,\,\,
\Delta = (\id \otimes \epsilon \otimes \id)\tau, \]

Then $\Delta$ is coassociative:
\[
(\id \otimes \Delta) \Delta (h) =
        \sum h^{(1)} \otimes \epsilon(h^{(2)})
        h^{(3)} \epsilon(h^{(4)}) \otimes h^{(5)}
= (\Delta \otimes \id) \Delta (h).
\]
Moreover, $\epsilon$ becomes a counit for coalgebra $(H,\Delta, \epsilon)$.
Indeed
\[\begin{array}{l}
 (\id \otimes \epsilon) \Delta(h)
= h^{(1)} \otimes \epsilon(\epsilon(h^{(2)})h^{(3)})\\
\,\,\,\,\,= h^{(1)} \otimes \epsilon(h^{(2)}) \epsilon(h^{(3)}) \\
\,\,\,\,\,=  h^{(1)} \otimes \epsilon(h^{(2)}h^{(3)}) \\
\,\,\,\,\,=  (\id \otimes \epsilon)(h^{(1)} \otimes h^{(2)}h^{(3)})
= (\id \otimes \epsilon) (h \otimes 1) = h \otimes \epsilon(1) \cong h,
\end{array}\]
\[\begin{array}{l}
 (\epsilon \otimes \id) \Delta(h)
= \epsilon(h^{(1)}\epsilon(h^{(2)})) \otimes h^{(3)}\\
\,\,\,\,\,= \epsilon(h^{(1)})\epsilon(h^{(2)})\otimes h^{(3)} \\
\,\,\,\,\,= \epsilon(h^{(1)}h^{(2)})\otimes h^{(3)}  \\
\,\,\,\,\,=  (\epsilon \otimes \id)(h^{(1)}h^{(2)}\otimes h^{(3)})
= (\id \otimes \epsilon) (1 \otimes h) = h \otimes \epsilon(1) \cong h.
\end{array}\]

As $\tau$ is an algebra map, we can easily see that
the map $\Delta$ is an algebra homomorphism too, so we have a bialgebra:
\[\begin{array}{lcl}
 \Delta(h) \Delta(g)& = &  \sum
(h^{(1)} \otimes \epsilon(h^{(2)}) h^{(3)})
(g^{(1)} \otimes \epsilon(g^{(2)}) g^{(3)}) \\
& = & \sum h^{(1)}g^{(1)} \otimes \epsilon(h^{(2)}) h^{(3)}
 \epsilon(g^{(2)}) g^{(3)} \\
& = & \sum h^{(1)}g^{(1)} \otimes \epsilon(g^{(2)}h^{(2)}) h^{(3)} g^{(3)} \\
& = & \sum (hg)^{(1)} \otimes \epsilon((hg)^{(2)}) (hg)^{(3)} \\
& = & \Delta(hg).\end{array}\]

We can also define the antipode
\[
Sh = \sum \epsilon(h^{(1)}) h^{(2)}  \epsilon(h^{(3)}).
\]
Indeed
\[\begin{array}{lcl}
\cdot(\id \otimes S)\Delta(h)& = & h^{(1)}S(\epsilon(h^{(2)})h^{(3)})
= h^{(1)}\epsilon(h^{(2)})S(h^{(3)})\\
&=& h^{(1)}\epsilon(h^{(2)})\epsilon(h^{(3)})h^{(4)}\epsilon(h^{(5)})\\
&=& [(\id \otimes \epsilon) (h^{(1)} \otimes h^{(2)}h^{(3)}) ]
h^{(4)}\epsilon(h^{(5)}) \\
&=&  [(\id \otimes \epsilon) (h^{(1)} \otimes 1) ]
h^{(2)}\epsilon(h^{(3)})\\
&=& h^{(1)} h^{(2)}\epsilon(h^{(3)})\\
&=& (\id \otimes \epsilon) (h^{(1)} h^{(2)}\otimes h^{(3)}) \\
&=& (\id \otimes \epsilon) (1\otimes h)\\
&=& \epsilon(h) 1_H
\end{array}\]
and similarily for $S$ at the left:
\[\begin{array}{lcl}
\cdot(S \otimes \id)\Delta(h)& = & S(h^{(1)}\epsilon(h^{(2)}))h^{(3)}
= Sh^{(1)}\epsilon(h^{(2)})h^{(3)}\\
&=& \epsilon(h^{(1)})h^{(2)}\epsilon(h^{(3)})\epsilon(h^{(4)})h^{(5)}\\
&=& \epsilon(h^{(1)})h^{(2)}
[(\epsilon \otimes \id) (h^{(3)} h^{(4)} \otimes h^{(5)}) ] \\
&=& \epsilon(h^{(1)})h^{(2)}[(\epsilon \otimes \id) (1 \otimes h^{(3)})]\\
&=& \epsilon(h^{(1)}) h^{(2)}h^{(3)}\\
&=& (\epsilon \otimes \id)(h^{(1)} \otimes h^{(2)} h^{(3)}) \\
&=& (\epsilon \otimes \id) (h\otimes 1)\\
&=& \epsilon(h) 1_H
\end{array}\]
Thus we have obtained a correspondence
from copointed quantum heaps into Hopf algebras where the underlying set is the same.
We leave to the reader to check that a map of underlying sets 
is a map of copointed quantum heaps
iff it is a map of Hopf algebras obtained via this correspondence. 
Thus the correspondence extends to a functor. 

\ppt {\bf Main theorem.} {\it The two functors constructed above are
mutually inverse isomorphisms of categories: copointed quantum heaps
$\Leftrightarrow$ Hopf algebras.
}

{\it Proof.} We need to show that the two functors are inverses. Underlying
algebra is identical, so we have to show that one composition
of the two functors does not change the coproduct $\Delta$ and
the other composition does not change the cooperation $\tau$.

Start with a copointed heap $(H,\mu,\eta,\tau,\epsilon)$. Then
\[\begin{array}{lcl}
\tau'(h) &=& (\id \otimes S \otimes \id)(\id \otimes \Delta) \Delta(h) \\
&=& (\id \otimes S \otimes \id) (\sum  h^{(1)} \otimes
\epsilon(h^{(2)})h^{(3)}\epsilon(h^{(4)})\otimes h^{(5)})\\
&=& \sum  h^{(1)} \otimes \epsilon(h^{(2)})
[\epsilon(h^{(3)})h^{(4)}\epsilon(h^{(5)})]\epsilon(h^{(6)})\otimes h^{(7)}\\
&=& \sum h^{(1)} \otimes \epsilon(h^{(2)}h^{(3)})h^{(4)}
\epsilon(h^{(5)}h^{(6)})\otimes h_{(7)}\\
&=&\sum [(\id \otimes \epsilon)(h^{(1)} \otimes h^{(2)}h^{(3)})]
h^{(4)} [(\epsilon\otimes \id)(h^{(2)}h^{(3)} \otimes h^{(4)})]\\
&=&\sum [(\id \otimes \epsilon)(h^{(1)}\otimes 1)]\cdot
h^{(2)}\cdot [(\epsilon\otimes \id)(1 \otimes h^{(3)})]\\
&=& \sum h^{(1)} \otimes h^{(2)} \otimes h^{(3)} = \tau(h).
\end{array}\]

Start with a Hopf algebra $(H,\mu,\eta,\Delta, \epsilon)$. Then
\[\begin{array}{lcl}
\Delta'(h)& = &\sum h^{(1)} \otimes \epsilon(h^{(2)}) h^{(3)} \\
&=& \sum h_{(1)} \otimes \epsilon(Sh_{(2)}) h_{(3)} \\
&=& \sum h_{(1)} \otimes \epsilon(h_{(2)}) h_{(3)} \\
&=& \sum h_{(1)} \otimes h_{(2)} = \Delta(h).
\end{array}\]

\section{The context and the work of Grunspan}

The present author has discovered the notion of quantum heap and
proved the main theorem of this article in Spring 2000, and this
entered then as Ch. 9 in his thesis~\cite{Skoda:thes}
on coset spaces for quantum groups.
The coset spaces were constructed there using coactions of
Hopf algebras and gluing using
noncommutative localizations~(\cite{Skoda:thes, banach,
Skoda:qbun, Skoda:cs, Skoda:nloc}).
This included the nonaffine generalization of torsors for Hopf
algebras.

Later, and independently, Grunspan~(\cite{grunspan})
considered an approach to torsors
via paralelogram approach and dualized this to the noncommutative setup.
He cites Kontsevich~(\cite{Konts:operads})
for using earlier this approach in commutative
affine case.
In other words, he studied certain coheap monoids in the symmetric monoidal
category of bimodules over a fixed commutative ``base'' algebra
over a ground field. When the base algebra is the ground
field, this is the same as our earlier introduced concept of quantum heap.
There are few differences however, in this case as well.
The first is minor, namely Grunspan introduced the axiomatics
with one additional axiom, but Schauenburg~\cite{schauenburg:torsors}
proved later that this axiom is superfluous.
The second difference is in the scope of work. Our main theorem
concerns the relation to Hopf algebras, namely the role of forgetting
and then reintroducing an algebra character. 
Grunspan overlooks this theorem, but proceeds with a study of
bitorsor picture, with a construction of a left and right
``automorphism'' quantum heaps, without a need to  
specify a character.

Schauenburg~(\cite{schauenburg:torsors}) proves
that Grunspan's torsors are essentially Hopf-Galois extensions.
My approach~(\cite{banach, Skoda:qbun}), to glue Hopf-Galois
extensions along coaction compatible noncommutative biflat localizations
is more general, in the sense that it gives a larger class
of objects which have the right to be called
noncommutative torsors.
One can extend this localization picture to Grunspan's formalism
by introducing the localizations compatible, with the ternary (co)operation $\tau$,
to develop a sort of gluing theory as well.
%We will address this issue elsewhere.

\section{Heapy categories}

Given a quantum heap $(H,\mu, \eta, \tau)$, the category of
representations of its underlying algebra inherits an additional
structure: a ternary product on objects. Namely, the tensor product of
three $H$-modules $A,B,C$ has also an $H$-action:
if $a\otimes b\otimes c \in A\otimes B\otimes C$
then $h (a \otimes b\otimes c) := h^{(1)}a\otimes h^{(2)}b\otimes h^{(3)}c$.
This triple tensor product of $H$-modules is the object part of
a categorical ternary product, which is a functor
${\mathcal C} \times {\mathcal C}^{\rm op}\times {\mathcal C}\to {\mathcal C}$.
We denote this product by $(A,B,C)\mapsto A\lozenge B\lozenge C$.
Then 
$$(Q_1\lozenge Q_2 \lozenge Q_3) \lozenge Q_4\lozenge Q_5
= Q_1\lozenge (Q_2 \lozenge Q_3 \lozenge Q_4)\lozenge Q_5
= Q_1\lozenge Q_2 \lozenge (Q_3 \lozenge Q_4\lozenge Q_5)$$
In fact, the equalities above are true only after natural
identifications, which are analogous to the MacLane's coherences
for monoidal categories. If it were a small category and
if we neglect the coherences, we see that essentially
the equality $(Q_1\lozenge Q_2 \lozenge Q_3) \lozenge Q_4\lozenge Q_5
= Q_1\lozenge Q_2 \lozenge (Q_3 \lozenge Q_4\lozenge Q_5)$
just says that this category is a cop in category of categories.
In general appropriate coherence isomorphisms are introduced, and we call 
such structures heapy categories. 
We will discuss their coherences properly in
the forthcoming work~\cite{Skoda:heapycat} as well as
their connections to a PRO for (co)heaps.
An important notion in this context is the notion of a unit for
a heapy category. It is an object $\bf 1$
such that the objects ${\bf 1} \lozenge Q \lozenge Q$,
${\bf 1}$ and $Q \lozenge Q\lozenge {\bf 1}$ are isomorphic
for each $Q$ (again, we should require coherent isomorphisms
with certain dinaturality properties).

A PRO is a strict monoidal category
whose object part is the set of natural numbers with
the addition as the tensor product. The addition of morphisms
does not need to be commutative though. The PRO for coheap monoids is
generated by a morphisms $t : 1 \to 3$, $e : 0\to 1$
and $d : 2\to 1$  which satisfy
the relations $(1 + t + 1)t = (2 + t)t = (t + 2)t$,
$(d + 1)t = e + 1 = 1 + e =(1 + d)t$, $d(d+1)=d(1+d)$, $d(1+e)=d(e+1)=1$. 
Clearly, usual heaps correspond to those strict monoidal
functors from its {\it opposite PRO} (=for heaps) to the cartesian category of sets,
for which $d^{\rm op}:1\to 2$ maps to the usual diagonal $a \mapsto (a,a)$
and $e^{\rm op}:1\to 0$ to the cancelling map $a\mapsto ()$.
Considering nonstrict monoidal functors to the category of categories, 
or instead, the techniques of~\cite{MLPare}, one can systematically 
introduce coherences in this setup.

As expected from the main theorem of this article,
each rigid monoidal category (having duals $Q\mapsto Q^*$)
gives rise to a unital heapy category via
$Q_1 \lozenge Q_2\lozenge Q_3 := (Q_1\otimes Q_2^*)\otimes Q_3$;
and the unit of the rigid monodial category may be equipped in
canonical way with coherences for the unit of a heapy category.
Conversely, a unital heapy category may be made monoidal
via $Q_1 \otimes Q_2 := Q_1 \lozenge {\bf 1} \lozenge Q_2$, again
with appropriate coherences~(\cite{Skoda:heapycat}).
In this way, a category of rigid monoidal categories is equivalent to
the category of unital heapy categories. It is interesting to
further study how much the rigid monoidal category depends on the
choice of unit; and which heapy categories have unit at all. The torsor
picture suggests that nonunital heapy categories may be of much
more interest than the nonunital monoidal categories are.

About the language: co in {\it cop} mimics co in coalgebra;
furthermore, in the dialect of Kent, according to OED,
cop is {\it a small heap of hay or straw}.

\end{document}